\documentstyle[12pt,leqno,amssymb,graphics]{amsart}
\newtheorem{thm}{Theorem}[section]
\newtheorem{lemma}[thm]{Lemma}
\newtheorem{prop}[thm]{Proposition}

\theoremstyle{definition}
\newtheorem{dfn}[thm]{Definition}

\theoremstyle{remark}
\newtheorem{remark}[thm]{Remark}

\begin{document}

\newcommand{\<}{\rightarrowtail}
\newcommand{\su}{\subseteq}
\newcommand{\pa}{\partial}

\newcommand{\x}{\times}
\newcommand{\I}{[0,1]}
\newcommand{\A}{{-Id_U}}
\newcommand{\B}{{-Id_{\E}}}

\newcommand{\bbox}{\hfill}

\newcommand{\R}{{\Bbb R}}
\newcommand{\Z}{{\Bbb Z}}
\newcommand{\E}{{{\Bbb R}^3}}
\newcommand{\C}{{{\Bbb Z}/2}}

\newcommand{\Q}{F\times [0,1]}
\newcommand{\s}{\Sigma F}

\newcommand{\p}{\pi_{1}}
\newcommand{\pp}{\pi_{2}}
\newcommand{\ppp}{\pi_{3}}

\newcommand{\U}{{
\left[
\begin{smallmatrix}
1 & 0 \\
0 & 1
\end{smallmatrix}
\right]}}
\newcommand{\V}{{
\left[
\begin{smallmatrix}
0 & 1 \\
1 & 0
\end{smallmatrix}
\right]}}
\newcommand{\T}{{
\left(
\begin{smallmatrix}
1 & 1 \\
0 & 1
\end{smallmatrix}
\right)}}

\newcommand{\e}{{
\left[
\begin{smallmatrix}
1 \\
0
\end{smallmatrix}
\right]}}

\newcommand{\ee}{{\left[
\begin{smallmatrix}
0 \\
1
\end{smallmatrix}
\right] }}

\newcommand{\eee}{{\left[
\begin{smallmatrix}
1 \\
1
\end{smallmatrix}
\right] }}

\newcommand{\0}{\times \{ 0 \}}
\newcommand{\1}{\times \{ 1 \}}
\newcommand{\h}{{1\over 2}}
\newcommand{\f}{{\hat{F}}}

\newcounter{numb}

\title{Quadruple Points of Regular Homotopies of Surfaces in 3-Manifolds}
\author{Tahl Nowik}
\address{Department of Mathematics, Columbia University, New York, 
NY 10027, USA.} 
\email{tahl@@math.columbia.edu}
\date{August 7, 1998}
\maketitle





\section{Introduction}

\begin{dfn}\label{d1}

Let $F$ be a (finite) system of closed surfaces and $M$ a 3-manifold.
A regular homotopy $H_t :F\to M,\ \ t\in\I$ will be called closed if 
$H_0 = H_1$. We will denote a closed \emph{generic} regular homotopy by CGRH.
The number mod 2 of quadruple points of a 
generic regular homotopy $H_t$ will be denoted 
by $q(H_t)$ ($\in \C$.)

\end{dfn}

Max and Banchoff in [MB] proved that any generic regular homotopy 
of $S^2$ in $\E$
which 
``turns $S^2$ inside out,'' has an odd number of quadruple points. 
The main point was showing that 
any CGRH 
of $S^2$ in $\E$
has an \emph{even} number of quadruple points. 
Goryunov 
in [G]
expresses this from the 
Vassiliev Invariants
point of view,  
as follows: Let $Imm (S^2 , \E)$ be the space of
all immersions of $S^2$ in $\E$, and let $\Delta \su Imm (S^2 , \E)$
be the subspace of all non-generic immersions.
Choose some generic immersion $f_0:S^2\to\E$ as a base immersion.
For any generic immersion $f:S^2\to\E$ let $Q(f)\in\C$ be 
defined as $q(H_t)$ where $H_t$ is some generic regular homotopy between
$f_0$ and $f$. There exists such an $H_t$ since 
$Imm (S^2 ,\E)$ is connected, and this is well defined since any
CGRH has $q=0$. Furthermore, 
since generic immersions do not have quadruple points,
$Q$ will be constant on each connected component of
$Imm (S^2 , \E) -\Delta$. 
[G] then raises the question whether such a $Q$ may be defined
for any surface in $\E$, that is, whether for any
CGRH of any surface in $\E$ the number of quadruple
points is 0 mod 2. 
$Q$ will then be specified by choosing 
one base immersion in each connected component
of $Imm(F,\E)$.

We begin this work with a short
alternative to the pictorial part of [MB].

We then answer the question of [G] to the affirmative in:

Theorem \protect\ref{t1}:
\emph{Let $F$ be a system of closed surfaces, and let $H_t:F\to\E$ be any CGRH.
Then $q(H_t)=0$.}

This phenomenon is not true in general for any
3-manifold in place of $\E$,
as we demonstrate in various examples in Section 4. 
However, we prove the following positive result:

Theorem \protect\ref{t3}:
\emph{
Let $M$ be an orientable irreducible 3-manifold with $\ppp (M) =0$.
Let $F$ be a system of closed orientable surfaces. 
If $H_t:F\to M$ is any CGRH in the regular homotopy class of an
embedding, then $q(H_t)=0$.}

Theorems \protect\ref{t1} and \protect\ref{t3} will both be proved
by reduction to the following more fundamental result, which will be 
proved first:

Theorem \protect\ref{t2}:
\emph{Let $M$ be any 3-manifold with $\pp (M) = \ppp (M) =0$.
Let $F$ be a system of closed surfaces and let $D\su F$ be a 
system of discs, one disc in each component of $F$.
If $H_t:F\to M$ is a CGRH that fixes $D$ then $q(H_t)=0$.}

In Section 5 we give an explicit formula for the above mentioned invariant
$Q$ for embeddings of a system of tori in $\E$. That is, 
if $F$ is a union of tori, then
for
every embedding $f:F\to\E$ we assign $Q(f)\in\C$
such that whenever $H_t$ is a generic regular homotopy between
two embeddings $f$ and $g$, then
$q(H_t) = Q(f) - Q(g)$.

\section{$S^2$ in $\E$}

Let $D\su S^2$ be a disc, 
and let $Imm_D (S^2 , \R^3)$
be the space of all immersions $f:S^2\to\R^3$ such that $f|_D$ is 
some chosen embedding.
$\p (Imm_D (S^2 , \R^3))$ is known to be $\Z$ (as will be explained bellow.) 
[MB] presented a specific CGRH which has 0 mod 2 quadruple points, and 
proved that it is the generator of $\p (Imm_D (S^2 , \R^3))$. 
The proof 
used a sequence of 19 drawings 
of intermediate stages of the homotopy.
Our following proposition is more general (and uses no pictures.)

Let $A:\E\to\E$ be the rotation $(x,y,z)\mapsto (-x,-y,z)$,
then $-A$ is the reflection $(x,y,z)\mapsto (x,y,-z)$.
Let $U=S^2 - int D$ be parametrized as the unit disc in the $xy$-plane,
and so the rotation $(x,y) \mapsto (-x,-y)$ acts on $U$.

\begin{prop}\label{t0}
Let $D$, $U$ and $A$ be as above. 
Let $f:S^2\to\R^3$ be
an immersion such that $f|_D$ is an embedding into the $xy$-plane,
and such that $f|_{U}(-x,-y) = A \circ f|_U (x,y)$. 
Let $H_t:S^2\to\E$ be a regular homotopy 
with $H_0 = f$, $H_1 = -A\circ f$ and which fixes $D$.

Define $H'_t : F \to \E$ 
by $H'_t (x,y) = -H_t (-x,-y)$ on $U$, and $H'_t$ is fixed on $D$.
Finally define $G_t = H_t * H'_t$ (where $*$ 
denotes concatenation from left to right.)

Then $G_t$ represents some odd power of the generator 
of $\p (Imm_D (S^2 , \R^3) , f)$.
\end{prop}

\begin{pf}

Given independent $v_1,v_2 \in \R^3$ let $K(v_1,v_2)\in SO_3$ be the 
unique matrix in $SO_3$ 
who's first two columns are obtained from $v_1,v_2$ by the Gram-Schmidt 
process. 

Given a representative $J_t$ of an element of $\p (Imm_D (S^2 , \R^3))$,
we define the following map $\bar{J}:S^3\to SO_3$,
which will be regarded as an element of $\ppp (SO_3)$. 
Parametrize $S^3$ as the quotient space of 
$U \x \I$ with identifications 
$((x,y),0) \sim ((x,y),1)$ for any $(x,y)\in U$, 
and $((x,y),t) \sim ((x,y),t')$ for any
$(x,y)\in \partial U$, $t,t'\in \I$.
And so $S^3$ contains one copy of $\pa U$ which we will still call $\pa U$.
Now let $\bar{J}:S^3\to SO_3$ be defined 
by $\bar{J}((x,y),t)= 
K( {\partial \over \partial x} J_t(x,y) , 
{\partial \over \partial y} J_t(x,y) )$.
This is well defined since $J_t$ fixes $D$ and so the 2-frame
$({\partial \over \partial x} J_t(x,y) , 
{\partial \over \partial y} J_t(x,y))$ is fixed with respect to $t$ 
on $\pa U = \pa D$.

Now, the map $J_t \mapsto \bar{J}$ induces an isomorphism
$\p (Imm_D (S^2 , \R^3)) \to \ppp(SO_3) = \Z$. 
(This follows from Smales' Theorem [S], and was also done 
in [MB]. Following the more general Smale-Hirsch Theorem,
we will perform the corresponding computation 
for surfaces and 3-manifolds
which may be non-orientable, in the proof of 
Theorem \protect\ref{t2} bellow.
We will then have no global 2-frame
on $F-D$ and no global 3-frame on $M$ to work with.)

Let 
$c: S^3 \to S^3$ be defined by 
$((x,y),t)\mapsto ( (-x,-y) \ , \ (t+\h) \ \hbox{mod} \ 1)$. 
This is an involution 
which is 
equivalent to the antipodal map.
Now $H'_t (x,y) = -H_t (-x,-y)$, so also
$H_t (x,y) = -H'_t (-x,-y)$.
Together this means that 
$G_t (x,y) = -G_{(t +\h) \ mod \ 1} (-x,-y)$.
Let $d J_t$ denote the pair
$({\partial \over \partial x} J_t , 
{\partial \over \partial y} J_t )$.
Applying this $d$ to both sides of the last identity gives 
$d G_t (x,y) 
= d G_{(t +\h) \ mod \ 1} (-x,-y) $.
Composing with $K$ gives:
$\bar{G} = \bar{G}\circ c$.

And so $\bar{G}$ induces a map $g$ from $S^3 / c$  $(= \R P^3)$ to 
$SO_3$  $ (= \R P^3)$.
We need to show that $\bar{G}$ represents an odd power of the generator of 
$\ppp (SO^3)$. This is equivalent to 
$g$ inducing an isomorphism 
$H^3(\R P^3, \Z / 2) \to H^3(\R P^3, \Z / 2)$. 
If $a \in H^1(\R P^3, \Z / 2)$ is the generator, then $a^3$ is the 
generator of $H^3(\R P^3, \Z / 2)$ and so it is enough to show that
$g$ induces an isomorphism on $H^1(\R P^3, \Z / 2)$, which 
is equivalent to it inducing an isomorphism on $\p (\R P^3)$. 
So this is what we will now show.

Let $X$ be the subgroup of $SO_3$ of matrices with third column $(0,0,1)$,
then $X$ is a circle.
$D$ remains fixed inside the $xy$ plane, and so the tangent planes 
of $S^2$ at points of $\partial U = \partial D$ are horizontal, so it 
follows that $\bar{G}(\partial U)\su X$ or $X'$, where $X'$ is the coset
of $X$ of all elements of $SO_3$ with third column $(0,0,-1)$. By
change of coordinates $(x,y)\mapsto (x,-y)$ on $U$ if needed
(this preserves the conditions of the theorem,) 
we may assume $\bar{G}(\partial U)\su X$.

We now claim that $\bar{G}|_{\pa U}:\pa U \to X$ is a map of degree 
$\pm 2$. (We are not orienting $\pa U$ and $X$ and so the sign is meaningless.)
We look at say $G_0=f$ for computing $\bar{G}$ on $\pa U$.
A loop going once around $X$ describes one full rotation of $\E$
around the $z$ axis.
And so we need to
check how many rotations 
does the
horizontal 2-frame $({\pa \over \pa x}f , {\pa \over \pa y}f)$  perform
with respect to the 2-frame $(1,0,0) , (0,1,0)$, 
when traveling once around $\pa U$.
But $(1,0,0) , (0,1,0)$ is globally defined on $f(D)$ and 
$({\pa \over \pa x}f , {\pa \over \pa y}f)$ is globally defined on 
$f(U)$ (in the sense of immersions.)
And so, since the Euler number of $TS^2$ ( $=\chi (S^2)$ ) is 2,
we must have 2 relative rotations. 

Since $c$ preserves $\pa U$ and $\bar{G}\circ c = \bar{G}$,
it follows that half of $\pa U$, (from some 
$p\in \pa U$ to $c (p)$,)
is mapped by $g$ with degree 1 onto $X$. But $X$ as a loop, is the 
generator of $\p (SO_3)$. And so the generator 
of $\p (S^3 / c)$ is mapped by $g$ to the generator
of $\p (SO_3)$, which is what we needed to show.

\end{pf}

The essential part of the Max-Banchoff theorem now follows:

\begin{thm}\label{mb1}
$q(J_t)=0$ for any CGRH $J_t:S^2\to \E$.
\end{thm}

\begin{pf}
There exists a generic $f$ satisfying the conditions of Proposition 
\protect\ref{t0}. 
(e.g. any embedded sphere of revolution which is parametrized as such.)
By [S] there exists a regular homotopy,
and thus also a generic regular homotopy $H_t$ from $f$ to $-A\circ f$
which fixes $D$.
The $G_t$ 
which is constructed 
from $H_t$
in Proposition \protect\ref{t0}, has exactly 
twice as many
quadruple points as $H_t$, and so $q(G_t)=0$. 
By Lemma 1 of [MB], $q$ is well defined on 
$\p (Imm_D (S^2 , \R^3) , f)$ (i.e. two CGRH's 
in the same class have the same $q$.)
Since $G_t$ represents an odd
power of the generator, and $q(G_t)=0$, we must also have $q=0$ for 
the generator, and thus $q=0$ for any element of 
$\p (Imm_D (S^2 , \R^3) , f)$.

Now let $J_t:S^2\to\E$ be any CGRH. By composing with a global isotopy 
of $\E$ we may assume $J_t|_D=f|_D$ for all $t$. Let $J'_t:S^2\to\E$
be a generic regular homotopy from $f$ to $J_0$ which fixes $D$.
Then $J'_t * J_t * J'_{-t}$ represents an element
of $\p (Imm_D (S^2 , \R^3) , f)$ and so 
$q(J_t) = q(J'_t * J_t * J'_{-t})=0$.

\end{pf}

To complete the picture, we follow [MB] from this point:
There are two isotopy classes of embeddings of $S^2$ in $\E$.
A regular homotopy from one 
of these isotopy classes
to the other, is called an \emph{eversion of the sphere}. 
By Theorem \protect\ref{mb1} any two generic 
eversions will have 
the same number mod 2 of quadruple points, 
so we need only to count this for one eversion. 
The Froissart-Morin eversion has exactly one 
quadruple point, and so we get the Max-Banchoff Theorem:

\begin{thm}\label{mb2}
Every generic eversion of the sphere has an odd number of quadruple points.
\end{thm}

\section{Surfaces in 3-Manifolds: Positive Results}

\begin{dfn}
Let $F$ be a system of closed surfaces, and $M$ a 3-manifold.
Two regular homotopies $H_t,G_t:F\to M$ will be called \emph{equivalent}
if as paths in the space $Imm(F,M)$,
they are homotopic relative their endpoints. 
(This is stronger than just having
the maps $H,G:\Q\to M$ defined by $H_t$ and $G_t$ be homotopic
relative $\pa (\Q)$.) 
\end{dfn}

\begin{lemma}\label{l1mb}
Let $F$ be a system of closed surfaces, and $M$ a 3-manifold.
If $H_t,G_t:F\to M$ are equivalent generic regular homotopies
then $q(H_t)=q(G_t)$.
\end{lemma}

\begin{pf}
The proof follows exactly as in Lemma 1 of [MB]. We only need to
emphasize, that for a generic regular homotopy $H_t$, by definition,
$H_0$ and $H_1$ are generic immersions, and so $H_0$ and $H_1$
do not have quadruple points.
\end{pf}

The following lemma is obvious:

\begin{lemma}\label{obv}
Let $F$ be a system of closed surfaces and $M$ a 3-manifold.
If $F^1,...,F^n$ are the connected components of $F$
then:

1. $Imm(F,M)=Imm(F^1,M)\x ... \x Imm(F^n,M)$.

2. Any regular homotopy $H_t:F\to M$ is equivalent to
a concatenation $H^1_t * ... * H^n_t$ where each $H^i_t$ fixes all 
components of $F$ except $F^i$.

\end{lemma}

As for the special case of $S^2$ in $\E$, we begin by assuming
that $H_t$ fixes discs:

\begin{thm}\label{t2}
Let $M$ be a 3-manifold with $\pp (M) = \ppp (M) =0$.
Let $F$ be a system of closed surfaces and let $D\su F$ be a 
system of discs, one disc in each component of $F$.
If $H_t:F\to M$ is a CGRH that fixes $D$ then $q(H_t)=0$.

\end{thm}

\begin{pf}

We assume first that $F$ is connected. The general case will
follow easily.

The Smale-Hirsch Theorem [H] states that $d:Imm(F,M) \to Mon(TF,TM)$ is a 
weak homotopy equivalence, 
where $Mon(TF,TM)$ is the space of bundle 
monomorphisms 
$TF\to TM$
(i.e. bundle maps 
which are a monomorphism on each fiber,)
and $d$ is the differential. 
From this it is easy to deduce, by means of [S] or [H], the following
relative version: Let $Imm_D(F,M)$ be the space of all immersions
$f:F\to M$ with $f|_D = H_0 |_D$ and let $Mon_D (TF,TM)$ be the space of all
bundle monomorphisms $f:TF\to TM$ with $f|_{TD} = dH_0 |_{TD}$,
then $d:Imm_D (F,M) \to Mon_D (TF,TM)$ is a weak homotopy equivalence.
In particular:
$d_* : \p (Imm_D (F,M), H_0) \to \p ( Mon_D (TF,TM), dH_0 )$ is an
isomorphism.

Our $H_t$ may be viewed as a representative of an 
element of $\p ( Imm_D (F,M) , H_0 )$. 
By Lemma \protect\ref{l1mb},
$q$ is well defined for elements of $\p ( Imm_D (F,M) , H_0 )$. 
We will show that $\p (Imm_D (F,M), H_0)$ ($= \p (Mon_D (TF,TM), dH_0)$)
is a cyclic group. 
We will then construct a CGRH $G_t$ with $q(G_t)=0$ and which is 
a representative of the generator of 
$\p (Imm_D (F,M) , H_0)$.
It will follow that $q=0$ for any generic representative of any element
of
$\p (Imm_D (F,M) , H_0)$, in particular for $H_t$.

Some definitions and notation:
Given a bundle map $f:TF \to TM$ let $\hat{f}:F\to M$
denote the map that it covers.
Denote by $Mon_{D,H_0} (TF,TM)$ the subspace 
of $Mon_{D} (TF,TM)$
of all $f\in Mon_{D} (TF,TM)$ with $\hat{f} = H_0$. 
Denote by $Map_D(F,M)$ the space of all maps $f:F\to M$ with
$f|_D = H_0 |_D$. 
For a space $X$ denote by $\Omega_x X$ the loop space
of $X$ based at $x\in X$. If $X$ is a space of maps $A\to B$
(as any of the above spaces,) and $a\in \Omega_f X$
then for $t\in\I$ the map $a(t)$ from $A$ to $B$ will be
denoted by $a_t:A\to B$.

It is shown in [HH] that any 3-manifold admits a connection, 
with parallel transport along any loop assigning either $Id$ or $-Id$.
Choose such a connection on $M$.
Given $a\in\Omega_{H_0} Map_D(F,M)$
it defines isomorphisms 
$A_{xt}: T_{a_0 (x)} M \to T_{a_t (x)} M$  
given
by 
parallel transport
along the path $s\mapsto a_s(x),\ s\in [0,t]$.
We will say $A_{xt}$ is associated to $a_t$.
We note that 
for any $x\in F$,
$A_{x1}$ is $Id$ and not $-Id$, 
since the loop
$t\mapsto a_t(x), \ t\in\I$ is nul-homotopic. 
This is so since it is homotopic to
$t\mapsto a_t(y), \ t\in\I$
with $y\in D$.

We now define a map 
$$\Phi:\Omega_{dH_0} Mon_D (TF,TM) \to \Omega_{H_0} Map_D (F,M)  
\x\Omega_{dH_0} Mon_{D,H_0} (TF,TM)$$ 
as follows: Given a loop $u\in \Omega_{dH_0} Mon_D (TF,TM)$, 
let $A_{xt}: T_{\widehat{u_0} (x)} M \to T_{\widehat{u_t} (x)} M$ be the 
continuous family of isomorphisms associated
to $\widehat{u_t}$.
Let $\widehat{\widehat{u_t}}:TF\to TM$ be defined by
$\widehat{\widehat{u_t}}|_{T_x F} = {A_{xt}}^{-1}\circ u_t|_{T_x F}$.
Now let $\Phi(u_t)=(\widehat{u_t} , \widehat{\widehat{u_t}})$.
$\Phi$ has the following inverse: given $(a,b)\in
\Omega_{H_0} Map_D (F,M)  \x\Omega_{dH_0} Mon_{D,H_0} (TF,TM)$,
take the $A_{xt}$ associated to $a_t$ and define ${\Phi}^{-1}(a,b)$
by $({\Phi}^{-1}(a,b))_t |_{T_x F}
={A_{xt}}\circ b_t|_{T_x F}$. And so $\Phi$ is a homeomorphism.

So $$\p (Mon_D (TF,TM), dH_0) = \pi_0(\Omega_{dH_0} Mon_D (TF,TM))$$
$$=\pi_0 
(\Omega_{H_0} Map_D (F,M)  
\x 
\Omega_{dH_0} Mon_{D,H_0} (TF,TM))$$
$$=\pi_0 (\Omega_{H_0} Map_D (F,M) )
\x 
\pi_0 (\Omega_{dH_0} Mon_{D,H_0} (TF,TM)). $$

We will first show that $\pi_0 (\Omega_{H_0} Map_D (F,M) )=0$.
An element in $\Omega_{H_0} Map_D (F,M) $ is a map
$h:\Q\to M$ such that $h(x,t)=H_0(x)$ whenever
$(x,t)\in F\0 \cup F\1 \cup D\x\I$. 
Let $SF$ be the quotient space of $\Q$ obtained by identifying 
$(x,0)\sim (x,1)$ for any $x\in F$ and $(x,t)\sim (x,t')$ for any
$x\in D$, $t, t'\in\I$. There is a natural inclusion $F=F\0 \su SF$.
A map $h:\Q\to M$ satisfying the above conditions 
is equivalent to a map
$h:SF\to M$ with $h(x)=H_0(x)$ for $x\in F$. We claim that any
two such maps are homotopic relative $F$, which means that
$\pi_0 (\Omega_{H_0} Map_D (F,M) )$ has just one element.
This is so, since there is a CW decomposition of $SF$ 
with $F$ being a subcomplex and 
$SF - F$ containing only 
2-cells and 3-cells, and since $\pp(M)=\ppp(M)=0$.

We now deal with
$\pi_0 (\Omega_{dH_0} Mon_{D,H_0} (TF,TM))$. 
For $x\in M$ let
$GL^+(T_x M)$ be the group of orientation preserving
automorphisms of $T_x M$ (i.e. automorphisms with positive determinant.) 
Let $GL^+TM$ be the bundle over $M$ who's fiber over $x$ is $GL^+(T_x M)$.
Let $Map_{D,H_0} (F, GL^+TM)$ be the space of all maps 
$f:F\to GL^+TM$
that lift $H_0$ and such that $f(x)=Id_{T_{H_0(x)} M}$ for every $x\in D$.

We define a map $\psi : Mon_{D,H_0} (TF,TM) \to Map_{D,H_0} (F, GL^+TM)$
as follows: 
Choose a metric on $M$.
Let $Mon_x$ denote the space of linear monomorphisms 
$T_x F\to  T_{H_0(x)} M$. 
We first
define a map $S_x :Mon_x\times Mon_x\to GL^+(T_{H_0(x)} M)$ as follows:
$S_x(\phi_1 , \phi_2)$ is defined to be the unique element in 
$GL^+(T_{H_0(x)} M)$ 
which extends
${\phi_2}\circ {\phi_1}^{-1} : {\phi_1}(T_x F)\to T_{H_0(x)} M$ 
and sends a unit vector normal to
$\phi_1(T_x F)$ into a unit vector normal to $\phi_2 (T_x F)$. 
Note that $S_x(\phi,\phi)=Id$ for any $\phi$.
Now, for $f\in Mon_{D,H_0} (TF,TM)$, define $\psi(f)(x) = 
S_x (dH_0|_{T_x F}, f|_{T_x F})$. 

We claim that $\psi$ is a homotopy equivalence. Its homotopy inverse 
$\psi' : Map_{D,H_0} (F, GL^+TM) \to  Mon_{D,H_0} (TF,TM) $
is given by
$\psi' (f) |_{T_x F} = f(x) \circ dH_0 |_{T_x F}$.
Indeed $\psi' \circ \psi = Id$. The homotopy 
$K:Map_{D,H_0} (F, GL^+TM)\x\I\to Map_{D,H_0} (F, GL^+TM)$
showing $\psi \circ \psi' \sim Id$
is given by $K(f,t)(x) = tf(x) + (1-t) \psi \circ \psi' (f) (x)$.
A priori, $tf(x) + (1-t) \psi \circ \psi' (f) (x)$
is just an endomorphism of $T_{H_0(x)} M$, but since 
$f(x)$ and $\psi \circ \psi' (f) (x)$ agree on 
the 2 dimensional subspace $dH_0(T_x F)$, and are both non-singular and
orientation preserving, then any convex combination of them will also be such.

So it is enough for us to consider $Map_{D,H_0} (F, GL^+TM)$.
But now we notice that $GL^+TM$ is a trivial bundle. Indeed, 
choose a base point $x_0\in M$. Our parallel transport
identifies any $T_x M$ with $T_{x_0} M$ only up to 
$\pm Id$, 
but since $-Id$ is in the center of $GL(T_{x_0} M)$,
the identification $GL^+(T_{x} M) \to GL^+(T_{x_0} M)$ is 
nevertheless well defined.
And so $Map_{D,H_0} (F, GL^+TM)$ is homeomorphic to the space 
$Map_D (F, GL^+(T_{x_0} M))$ of all maps $f:F\to GL^+(T_{x_0} M))$ 
with $f(x)=Id$ for every $x\in D$. Denote $GL^+=GL^+(T_{x_0} M)$.
We were interested in 
$\pi_0 (\Omega_{dH_0} Mon_{D,H_0} (TF,TM))$. This corresponds now to 
$\pi_0 (\Omega_{Id} Map_D (F, GL^+) )$ 
where $Id$ here is the constant map taking 
$F$ into $Id\in GL^+$. This in turn is the same as the set of 
based homotopy classes
$[(\Sigma F,*), (GL^+,Id)]$ where $\Sigma F$ is the quotient space
of $\Q$ with $F\0 \cup F\1 \cup D\x\I$ identified into one point 
$*$.

What we obtained so far, is an isomorphism
$\p (Imm_D (F,M), H_0) \to [(\Sigma F,*), (GL^+,Id)]$. 
But ${GL_n}^+$ is homotopy equivalent to $SO_n$, so we
are finally interested in $[(\Sigma F,*), (SO_3,Id)]$. 
(The group structure on $[(\Sigma F,*), (SO_3,Id)]$ 
is given by concatenation along the $\I$ variable of $\Sigma F$.)

Now, there is a CW decomposition of $\Sigma F$ with one 0-cell,
no 1-cells, some number of 2-cells and one 3-cell. 
Let $S^3$ be modelled as the quotient space of $\Sigma F$ with the 
2-skeleton of $\Sigma F$ identified into one point $*$, and let
$e:\Sigma F \to S^3$ be the quotient map.
Since 
$\Sigma F$ has no 1-cells
and $\pi_2(SO_3)=0$, any map $\Sigma F \to SO_3$ is homotopic
to a map which factors through $S^3$.
In other words, the map 
$e^* : [(S^3,*), (SO_3,Id)] \to [(\Sigma F,*), (SO_3,Id)]$ 
is an epimorphism.  
But $[(S^3,*), (SO_3,Id)] = \ppp(SO_3) = \Z$. and so 
$\p (Imm_D (F,M), H_0) = [(\Sigma F,*), (SO_3,Id)]$ 
is cyclic. 

We construct the following CGRH $G_t:F\to M$ which fixes $D$ and has
$G_0=H_0$.
Let $B\su M$ be a small ball such that $B\cap H_0(F)$ is an embedded disc, 
disjoint from $H_0(D)$.
Parametrize $B$ such that $H_0(F)\cap B$ looks in $B$ like
$(s-d) \cup t$ where $s$ is a sphere in $B$, $t$ is a very thin tube
connecting 
$s$ to $\pa B$, and $d\su s$ is the tiny disc deleted from $s$ 
in order to glue $t$.
Let $i:s\to B$ denote the inclusion map.
Let $G'_t$ be a generating CGRH of $\p (Imm_d(s,B) , i)$.
By changing $G'_t$ slightly, we may also assume that
no quadruple point occurs in $d$, and no triple point passes $t$. 
Now, identify $u=s-d$ with its preimage in $F$, and define $G_t:F\to M$ 
as $G'_t$ on $u$ and as fixed on $F-u$.
The conditions on $G'_t$ guarantee that the number of quadruple points of 
$G'_t:s\to B$ and $G_t:F\to M$ is the same, and so 
$q(G_t)=q(G'_t)=0$ by Theorem \protect\ref{mb1}.
Our proof (for the case $F$ connected)
will be complete if we show $G_t$ is a generator of $\p (Imm_D
(F,M) , H_0)$.

The above mentioned CW decomposition of $\Sigma F$ comes from the product
structure on
$F\x\I$ and so there is no problem choosing $B$ such that 
$int\ u\x\I\su\Sigma F$ 
will be contained in the open 3-cell of $\Sigma F$,
and so our model of $S^3$ will contain a copy of $int\ u\x\I$.
As done for $F$, $\Sigma s$ will denote the quotient space of $s\x\I$ with
$s\0 \cup s\1 \cup d\x\I$ identified into one point $*$.
$\Sigma s$ also contains a copy of $int\ u\x\I$.
Define $f:S^3\to \Sigma s$ ($\cong S^3$) by 
$f(x)=x$ for all $x\in int\ u\x\I$, and 
$f(S^3 - int\ u\x\I)=*$.
$f$ is clearly of degree 1, and so $f^* : [(\Sigma s,*) , (SO_3,Id)] \to
[(S^3,*),(SO_3,Id)]$
is an isomorphism.  
And so  
$(f\circ e)^* : [(\Sigma s,*) , (SO_3,Id)] \to [(\Sigma F,*),(SO_3,Id)]$
is an epimorphism. Denote by $\bar{J}$ the map $\Sigma F \to SO_3$
which we have attached to a CGRH $J_t$.
For $\bar{G'}$ 
to be defined we must choose a connection and metric on $B$. If
we 
use the restrictions to $B$ of the connection and metric chosen for $M$,
then
we get $(f\circ e)^*(\bar{G'}) = \bar{G}$.
And so $G'_t$ generating
$\p (Imm_d (s,B) , i)$ implies $\bar{G'}$ generates $[(\Sigma s ,*) ,
(SO_3,Id)]$ 
implies $\bar{G}$ generates $[(\Sigma F ,*) , (SO_3,Id)]$ implies finally
that $G_t$ generates $\p (Imm_D (F,M) , H_0)$.

We now deal with the general case where $F$ may be non-connected.
Let $F^1,...,F^n$ be the components of $F$ and $D^i\su F^i$ the components
of $D$. As in Lemma \protect\ref{obv},
$Imm_D(F,M) = Imm_{D^1}(F^1,M)\x ... \x Imm_{D^n}(F^n,M)$
and so 
$$\p(Imm_D(F,M),H_0) = \p(Imm_{D^1}(F^1,M),H_0|_{F^1})
\x ... \x \p(Imm_{D^n}(F^n,M),H_0|_{F^n}).$$ So
$\p(Imm_D(F,M),H_0)$ is generated by the $n$ generators
$G^i_t\in \p(Imm_{D^i}(F^i,M), H_0|_{F^i})$. But 
by our proof,
we may choose the generator 
$G^i_t$ to move $F^i$ only
in a small ball $B^i$ which does not intersect any other component.
And so the number of quadruple points of $G^i_t$ when thought of 
as an element of $\p (Imm_{D^i}(F^i,M),H_0|_{F^i})$ or $\p (Imm_D(F,M),H_0)$ 
is the same. Since $q(G^i_t)=0$, and $G^i_t$ are generators, the theorem
follows.
\end{pf}

Our aim from now on will be, to be able to reduce a general CGRH, 
to a CGRH
that fixes such a system of discs.

\begin{lemma}\label{ai}
Let $M$ be any 3-manifold. Let $\delta\su M\x M$ be the diagonal, i.e.
$\delta = \{(x,x) \ : \ x\in M \}$. There exists an open neighborhood
$\delta \su U \su M\x M$ and a continuous map 
$\Phi: U\to Diff(M)$ where $Diff(M)$ is the space of 
self diffeomorphisms of $M$, such that:

1. $\Phi (x,x) = Id_M$ for any $x\in M$.

2. $\Phi (x,y)(x) = y$ for any $(x,y) \in U$.

\end{lemma}

\begin{pf}
Let $B_r\su \E$ denote the ball of radius $r$ (about the origin.)
Choose an isotopy $f_r:\E \to \E$, $r\in [0,\infty )$ 
with the following properties: (1) $f_r$ is the identity outside $B_{2r}$
(2) $f_r(0,0,0) = (0,0,r)$
(3) For any $A\in O_3$ which fixes the $z$-axis, $f_r \circ A = A \circ f_r$.

Choose a metric on $M$. 
Let $d(x,y)$ denote the distance in $M$ and
let $B_r(x)$ denote the ball of radius $r$ about $x\in M$.
Let $\delta \su U \su M\x M$ be a thin neighborhood such that
for any $(x,y)\in U$, the ball in $T_x M$ of radius $2d(x,y)$ is
embedded by the exponential map.

For $(x,y)\in U$ let $r=d(x,y)$, and define $\Phi (x,y)$ as follows: 
$\Phi (x,y)$ will be the identity outside $B_{2r}(x)$. In $B_{2r}(x)$ take 
normal coordinates 
(i.e. an orthonormal set of coordinates on $T_x M$ which is 
induced onto $B_{2r}(x)$ via the exponential map,)
with the additional requirement that the $z$ axis runs from $x$ to $y$, i.e.
that in these coordinates $x=(0,0,0)$ and $y=(0,0,r)$.
We let $\Phi (x,y)$ act on $B_{2r}(x)$ 
as $f_r$ with respect to these coordinates.
Due to the symmetry of $f_r$, this definition does not depend on the 
freedom we still have in the choice of the normal coordinates. 
It is easy to verify that $\Phi:U\to Diff(M)$ is indeed continuous. 
\end{pf}

\begin{lemma}\label{nh}
Let $M$ be any 3-manifold and $F$ a (connected) closed surface.
Let $H_t:F\to M$ be a regular homotopy. Let $p\in F$ 
and let $h:\I\to M$ be defined by $h(t)=H_t(p)$.

If $g:\I\to M$  is homotopic to $h$ relative end points, then there 
exists a regular homotopy
$G_t:F\to M$ which is equivalent to $H_t$ and with $G_t(p)=g(t)$.

In particular, if $H_t$ is closed, and $h(t)$ is nul-homotopic, then
$H_t$ is equivalent to a CGRH which fixes $p$.

\end{lemma}

\begin{pf}
Let $k:\I\x\I\to M$ satisfy 
$k(t,0)=h(t)$, $k(t,1)=g(t)$, $k(0,s)=h(0)$, $k(1,s)=h(1)$.
Let $U$ and $\Phi (x,y)$ be as defined in
Lemma \protect\ref{ai}.
There is a partition $0=s_0 < s_1 < ... < s_n = 1$ such that
$( k(t,s_{i-1}) , k(t,u) ) \in U$ for any $i$, $t$ 
and $s_{i-1} \leq u \leq s_i$. 

Let $I^i:\I\x [s_{i-1},s_i]\to Diff(M)$
be defined as $I^i(t,u)=\Phi (k(t,s_{i-1}),k(t,u) )$.
Extend $I^i$ to $\I\x\I$ by letting $I^i(t,u)=Id_M$ for $0\leq u \leq s_{i-1}$
and $I^i(t,u)= I^i(t,s_i)$ for $s_i \leq u \leq 1$.
Now let $K(t,u)=I^n(t,u) \circ ... \circ I^1(t,u) \circ H_t$,
then $K(t,0)=H_t$ and
$G_t=K(t,1)$ is the required equivalent regular homotopy.

\end{pf}

\begin{lemma}\label{fp}
Let $M$ be any 3-manifold, $F$ a (connected) closed surface and $p\in F$.
If $H_t$ is a CGRH that fixes $p$, and such that $p$ is not in 
the multiplicity set of $H_0$,
then there exists a 
generic regular homotopy 
$G_t$ and discs $D,D'$ with $p\in D \su intD' \su F$ such that:

1. $G_0 = H_0$.

2. $G_t$ fixes $D$.

3. Either $G_1 = H_1$ ($=H_0$) or $G_1 = H_1 \circ d$ where $d:F\to F$ 
is a Dehn twist performed on the annulus $A = D'-D$.

4. $q(G_t) = q(H_t)$.

\end{lemma}

\begin{pf}

Let $p'=H_t(p)$, $t\in [0,1]$.
By slight flattening out and straightening out we may assume that
there are coordinates on a neighborhood 
$p'\in U \su M$
(with the origin corresponding to $p'$,) 
and a small disc
$p\in D' \su F$ such that 
$H_t(D')\su U$, $t\in\I$ and 
such that 
with respect to these coordinates on $U$ the following holds:
(1) $H_0(D')$ is an actual flat round disc (2)
$H_t |_{D'} = I_t \circ H_0|_{D'}$
where $I_t \in SO_3$, $I_0 = I_1 = Id$. 

We continue to work in the chosen coordinates, and use the natural norm $||x||$
of the coordinates.
We may assume the radius of 
$H_0(D')$ is 2, and let $D\su D'$ be the disc with $H_0(D)$ having radius 1.
Let $B_r \su U$ denote the ball of radius $r$, 
then by our assumption on $p$ we may also assume $B_2\cap H_0(F) = H_0(D')$.
Now let $J_t:M \to M$ be 
the following isotopy:
On $B_1$, $J_t={I_t}^{-1}$. For $1 \leq ||x|| \leq 2$, 
$J_t(x) = ({I_{(2-||x||)t}})^{-1}(x)$. 
For $x\in M-B_2$, $J_t(x)=x$ for all $t$. 
Let $H'_t = J_t\circ H_t$. Then $H'_t$ satisfies
1,2 and 4 of the lemma, and $H'_1 = H_1$ ($=H_0$) on $F-A$.
For $x\in A$, $H'_1(x) = {(I_{2-||H_1(x)||})}^{-1} (H_1 (x))$.
The map $k:[1,2]\to SO_3$ defined by $k(s)={(I_{2-s})}^{-1}$
has $k(1)=k(2)=Id$, and so there is a homotopy $K: [1,2] \x [0,1] \to SO_3$
with $K(s,0) = k(s)$, $K(1,u)=K(2,u)=Id$ and $K(s,1)$ is either
the constant loop on $Id$, or the loop 
that describes one full
rotation about the axis perpendicular to $H_1(D')$.
(The loops here have domain $[1,2]$.)
We now define $H''_t :F\to M$ as follows:
For $x\in F-A$, $H''_t(x)=H_1(x)$. For $x\in A$,
$H''_t(x)=K(||H_1(x)|| , t) (H_1 (x))$.
Finally, our desired $G_t$ is $H'_t * H''_t$.

\end{pf}

We combine Lemmas \protect\ref{l1mb}, \protect\ref{obv}, \protect\ref{nh}
and the proof of \protect\ref{fp} to get:

\begin{prop}\label{p0}
Let $M$ be any 3-manifold and
let $F=F^1\cup ... \cup F^n$ 
be a system of closed surfaces.
Let $H_t:F\to M$ be a CGRH
and let $p^i\in  F^i$.

If for each $i$, the loop $t\mapsto H_t(p^i)$ is nul-homotopic,
then there is a generic regular homotopy
$G_t:F\to M$ and discs $D^i \su int{D'}^i \su F^i$
such that:

1. $G_0 = H_0$.

2. $G_t$ fixes $\bigcup_i D^i$.

3. For each $i$,
either $G_1|_{F^i} = H_1|_{F^i}$ ($=H_0|_{F^i}$) 
or $G_1|_{F^i} = H_1|_{F^i} \circ d$ 
where $d:F^i\to F^i$ 
is a Dehn twist performed on the annulus $A^i = {D'}^i-D^i$.

4. $q(G_t) = q(H_t)$.

\end{prop}

\begin{pf}
If $t\mapsto H_t(p^i)$ is nul-homotopic, then the same is true for
any other $p\in F^i$, and so we may assume $p^i$ is not in the
multiplicity set of $H_0$.
By Lemmas \protect\ref{l1mb}, \protect\ref{obv} and \protect\ref{nh}
we may assume that $H_t=H^1_t * ... * H^n_t$ with each
$H^i_t$ fixing $p^i$ and fixing all $F^j$, $j\neq i$.
Now replace each $H^i_t$ with $G^i_t$ 
which fixes all $F^j$, $j\neq i$, and such that $G^i_t|_{F^i}$ is the 
regular homotopy constructed from $H^i_t|_{F^i}$
as in the proof of Lemma
\protect\ref{fp},
and making sure that each ${B_2}^i$ of that proof
satisfies ${B_2}^i \cap H_0(F) = H_0({D'}^i)$
(not just the intersection with $H_0(F^i)$.)

When it is the turn of $F^i$ to move, then $F^j$ for $j<i$ have already
performed their movement, but
notice that whether a Dehn twist appeared or
not, $G^j_1(F^j)=H^j_1(F^j)$. 
And so since ${H'}^i_t$ differs from
$H^i_t$ only in $B_2^i$, and 
${H''}^i_t$ moves only in $B_2^i$,
we have $q(G^i_t)=q(H^i_t)$. And so finally $q(G_t)=q(H_t)$.

\end{pf}

We are now ready to prove:

\begin{thm}\label{t1}
Let $F$ be a system of closed surfaces, and let $H_t:F\to\E$ be any CGRH.
Then $q(H_t)=0$.
\end{thm}

\begin{pf} 
Replace $H_t$ with 
the $G_t$ of Proposition \protect\ref{p0}. 
($t\mapsto H_t(p^i)$ are of course nul-homotopic.)
For brevity, we will denote both
$F^i$ itself, and its immersed image, by $F^i$. 
Denote $\hat{F}^i=\bigcup_{j\neq i} F^j$.
We continue $G_t$ to get a CGRH as follows. Each $F^i$ for which
$G_1|_{F^i}=H_1|_{F^i}$ 
will remain fixed. If for some $F^i$ a Dehn twist appeared,
then we undo this Dehn twist by rotating $F^i - {D'}^i$ by a rigid rotation 
about the line $l_i$ perpendicular 
to $F^i$ at $p^i$, and while keeping $D_i$ and $\hat{F}^i$ 
fixed. (We may assume that $l_i$ is generic with respect to 
$F$ in the sense that it intersects it generically, and that this 
rigid rotation of $F^i$ while keeping $\hat{F}^i$ fixed, is a generic
regular homotopy of $F$.) We do this one by one, to each one of the components 
$F^i$ for which a Dehn twist appeared. We claim that each such rigid rotation
contributes 0 mod 2 quadruple points. 
Indeed, by our genericity assumption,
such a quadruple point may occur in one of three ways: A triple point,
double line or sheet, of $F^i$ crosses respectively a sheet,
double line or triple point of $\hat{F}^i$.
We will show that each of the three types separately contributes 0 mod 2 
quadruple points: 
(1) A triple point of $F^i$ traces a circle during 
this rigid rotation. Since $\E$ is contractible, the $\C$ intersection of
any $\C$ 1-cycle (this circle) with any $\C$ 2-cycle ($\hat{F}^i$)
is 0. (2) Each double line of $F^i$ is an 
immersed loop, and so it traces a (singular) torus during the rigid rotation.
Again, the $\C$ intersection of this torus, with each immersed double loop
of $\hat{F}^i$ is 0. (3) This is symmetric to case 1. (One may think of
$\hat{F}^i$ as rotating around $l_i$ and $F^i$ remaining fixed.)

And so we have managed to complete $G_t$ to a CGRH without changing 
$q$. But this CGRH fixes $D=\bigcup_i D^i$. And so by Theorem 
\protect\ref{t2} we are done.

\end{pf}

We proceed towards Theorem \protect\ref{t3}.

\begin{figure}[h]
\ \ \ \ \ \ \ \ \ \ \ \ \ \ \ \ \ \ \ \ \ \ \ \ \   
\scalebox{.6}{\includegraphics{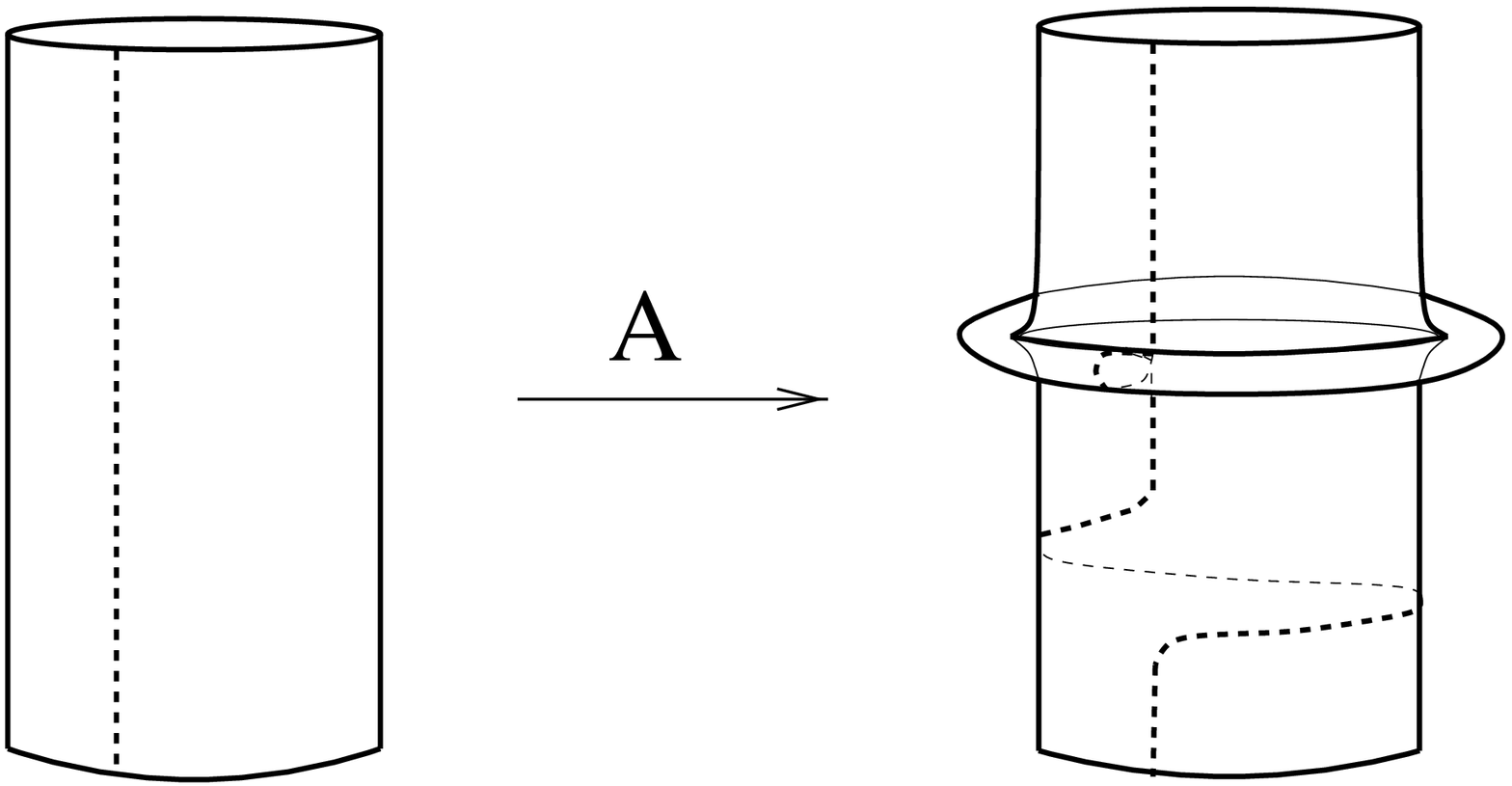}}
\caption{Move A}
\end{figure}

\begin{figure}[h]
\ \ \ \ \ \ \ \ \ \ \ \ \ \ 
\scalebox{.6}{\includegraphics{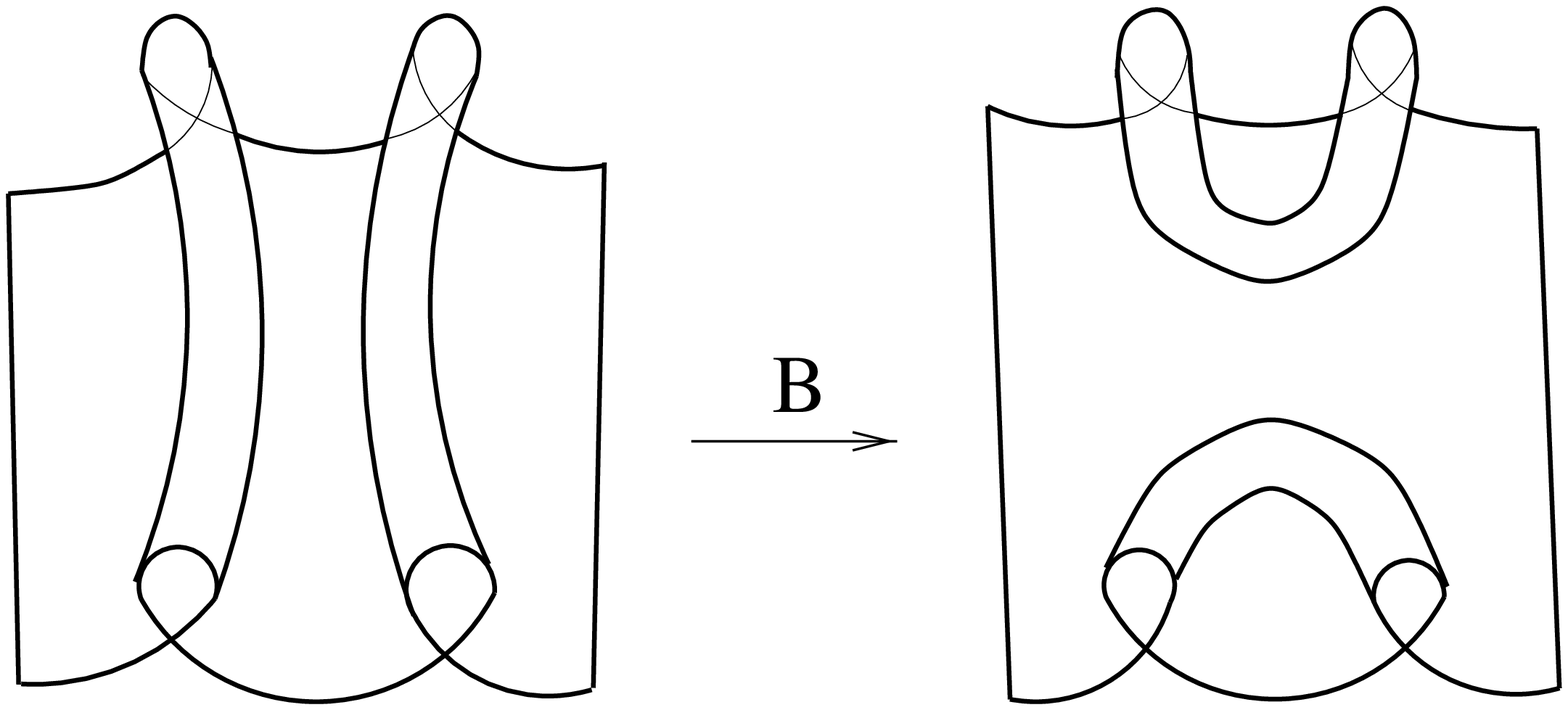}}
\caption{Move B}
\end{figure}

Figure 1 (resp. 2) describes a regular homotopy of an 
annulus (resp. disc) which is properly immersed in a ball (the
ball is not drawn.) The regular homotopy is assumed to fix a neighborhood 
of the boundary of the annulus (resp. disc.) This regular homotopy
will be called move $A$ (resp. move $B$.) Move $A$ 
begins with an embedding, and adds one circle of intersection, and 
one Dehn twist. Move $B$ begins with two arcs of intersection, and
replaces them with two different arcs of intersection.

We do not specify the regular homotopies themselves, we only specify
the initial and final immersions. It is easy to see that the 
desired regular homotopies exist, by the Smale-Hirsch Theorem: For 
move $B$ there is nothing to check since any two immersions of a disc $D$
in a ball $B$ which coincide on a neighborhood of $\pa D$, are 
regularly homotopic fixing a neighborhood of $\pa D$. For move $A$, let $K$
denote the annulus.
One needs to check that the map $K \to SO_3$ 
associated to the initial immersion, is homotopic relative $\pa K$
to the one associated to the final immersion. Since $\pp(SO_3)=0$,
we need to check this only for some spanning arc, say the dotted one,
in Figure 1.

Moves $A$ and $B$ will be applied as local moves to surfaces in 3-manifolds,
i.e. performing move $A$ or $B$ on a surface means that this regular homotopy 
is performed in a small ball in $M$, and the rest of $F$ remains fixed.

If $c$ is the essential circle of the annulus $K$, then we will say that
``move $A$ was applied to the circle $c$.'' 
If $f:F\to M$ is an immersion, then by definition, 
move $A$ may be applied to a 
circle $c \su F$ iff there is an (embedded) disc $E\su M$ 
with $E\cap f(F) =\pa E =f(c)$. The move is then performed 
in a thin neighborhood
of $E$.
In particular, move $A$ may always be performed
on a circle
$c\su F$ which bounds a disc $D$ in $F$, with $D$ 
containing no multiple points of 
$f$. 
In this case we may undo the Dehn twist simply by rotating $D$. 
The torus that is added
to $f(F)$ by an $A$ move will be called ``the \emph{ring} formed by the 
$A$ move.''  There are several choices to be made when applying an $A$ move
to a circle $c$. One may choose on what side of the ring will the Dehn twist
appear, what the orientation of the Dehn twist will be, and to what side of
$F$ will the ring be facing. (In Figure 1 the ring happens to be facing to
the side away from the disc $E$.)

\begin{lemma}\label{l1}
$q(A)=q(B)=1$ i.e. if $K$ denotes an annulus, $D$ denotes a disc, 
and $B$ denotes a ball, then for any generic
regular homotopy $A_t:K\to B$ that realizes the $A$ move, $q(A_t)=1$,
and for any generic regular homotopy $B_t:D\to B$
that realizes the $B$ move, $q(B_t)=1$.
\end{lemma}

\begin{pf}
For $A$: Take an embedding $f:S^2\to\E$. Take a ball $B\su \E$
such that $K = f(S^2)\cap B$ is a standard annulus in $B$. Perform $A_t$ on 
$K$ in $B$, while keeping the rest of $f(S^2)$ fixed, and
such that the ring obtained from this $A$ move will
face the outside of $f(S^2)$. 
This regular homotopy of $S^2$ may be continued
to complete an eversion of $S^2$ in $\E$
with no additional quadruple points. 
And so by Theorem \protect\ref{mb2} we must have 
$q(A_t)=1$.

For $B$: We describe the above sort of argument in shorter form:
Start again with an embedding $S^2\to\E$. 
Perform two $A$ moves at two remote regions of $S^2$, with both rings 
facing the outside of $S^2$.
Perform 
one $B$ move to merge the two rings into one. 
We may continue as before with no additional quadruple points 
to get an eversion. And so $2q(A)+q(B)=1$ and so $q(B)=1$.

\end{pf}

\begin{lemma}\label{l2}
Let $F$ be a (connected) closed surface.
Let $i:F=F\0 \to F\x [-1,1]$ be the inclusion. Let $D,D'$ be two discs 
$D\su intD' \su F$. Then, there exists a generic regular homotopy 
$H_t:F\to F\x [-1,1]$ which
satisfies:

1. $H_0=i$ 

2. $H_1=i\circ d$ where $d:F\to F$ is a Dehn
twist performed on the annulus $D'-D$. (We may choose the orientation 
of the Dehn twist.)

3. $H_t$ fixes $D$.

4. $q(H_t)=\chi (F) \ mod \ 2$. 

\end{lemma}

\begin{pf}
Let $E\su F\x [-1,1]$ be a disc with $\pa E = F\cap E$ an essential loop in 
$D'-D$. Perform move $A$ 
in a thin neighborhood of $E$ disjoint from $D$.
Now choose a Morse function on $F-int D'$,
$\pa D'$ being the minimal level curve. Let the ring created by 
the $A$ move, 
move along with the increasing level curves. Whenever there is a singularity
of the Morse function, namely a local minimum, local maximum or 
saddle point, then we will perform move $A$, ${A}^{-1}$ or $B$ respectively
in order for the rings to continue
following up with the level curve. 
If in all $A$ moves, we create the ring 
on the same side of $F$, then the $B$ moves will always be possible.
This we do until we pass the 
maximum of the Morse function, and so we remain with no
rings at all. Only the Dehn twist of the initial $A$ move remains,
since all others may be resolved by rotating the relevant disc.
(The relevant disc for the initial $A$ move is $D$, which we 
are keeping fixed.)
So the count of singularities will match the count of $A,A^{-1},B$ moves.
Together with our initial $A$ move to match the disc $D'$, we have
by Lemma \protect\ref{l1}, and by the relationship between singularities of
Morse functions and Euler characteristic:
$q(H_t)=\chi(F) \  mod \ 2$. 
\end{pf}

\begin{prop}\label{p1}
Let $M$ be a  3-manifold with $\pp (M) = \ppp (M) = 0$.  
Let $F=F^1 \cup ... \cup F^n$ be 
a system of closed
surfaces each having even Euler characteristic. Let $H_t:F\to M$ be a CGRH
in the regular homotopy class of a two sided embedding. 
Let $p^i\in F^i$. 
If for each $i$, the loop $t\mapsto H_t(p^i)$, $t\in\I$ is nul-homotopic,
then $q(H_t)=0$
\end{prop}

\begin{pf}
Let $J_t:F\to M$ be a generic regular homotopy from an embedding to $H_0$.
By passing to $J_t * H_t * J_{-t}$ 
we may assume $H_0$ is a (two sided) embedding.
Replace $H_t$ by the $G_t$ of Proposition \protect\ref{p0}.
For each $i$ for which there appeared a Dehn twist, 
we cancel the Dehn twist by concatenating with the regular homotopy of
Lemma \protect\ref{l2}, performed in a thin neighborhood of $H_0(F^i)$.
By our assumption on the Euler characteristic
of the $F^i$s, this does not change $q$.
And so we replaced $H_t$ with a CGRH that fixes $\bigcup_i D^i$
and has the same $q$.
We are done by Theorem \protect\ref{t2}.

\end{pf}

\begin{lemma}\label{cd}
Let $M$ be a 3-manifold with $\pp (M) = 0$. Let $F\su M$ be an embedded
system of closed surfaces, and let $D\su M$ be a compressing disc for 
$F$ i.e. $D\cap F = \pa D$. 
Let $F'\su M$ be the system of surfaces 
obtained from $F$ by compressing along $D$, in a thin neighborhood
of $D$. 
Let $i:F\to M$, $i':F'\to M$
denote the inclusion maps. 

If $i':F'\to M$ has the property that any CGRH $G_t:F'\to M$ with
$G_0=i'$ has $q(G_t)=0$, then $i:F\to M$ has the same property.
\end{lemma}

\begin{pf}
Let $H_t$ be a CGRH with $H_0 = i$. 
By slightly relocating $D$, we may assume $H_t$ is generic
with respect to $\pa D$ (considered here as a loop in $F$.) 
Let $\hat{F}$ denote the 2-complex
$F\cup D$. By means of [S] or [H], $H_t$ may be extended to a generic
regular homotopy $H'_t:\hat{F}\to M$, 
with $H'_0$ the inclusion and
such that $H'_1$ coincides with $H'_0$ on $F\cup N(\pa D)$
where $N(\pa D)$ denotes a neighborhood of $\pa D$ in $D$. 
Since $\pp (M) = 0$ we may further 
continue $H'_t$ (as a generic regular homotopy) 
until $H'_0=H'_1$ on the whole of $\hat{F}$. This regular homotopy
of $\hat{F}$ induces a CGRH $G_t:F'\to M$ with $G_0=i'$
in an obvious way.
By assumption, $q(G_t)=0$. 
$H'_t$ being a \emph{generic} regular homotopy of $\hat{F}$ means 
in particular that
any quadruple point of $H'_t$ occurs away from $\pa D$. 
And so there are two types of quadruple points of $H'_t$: 
(1) The quadruple points of $H_t$. (2) Quadruple points which involve at
least one sheet coming from $D$.
We assume the compression along $D$ was performed in a very thin 
neighborhood of $D$ and that $G_t$ follows $H'_t$ very closely, and so
$G_t$ will inherit the quadruple points of $H'_t$ in the following way:
Every quadruple point of type (1) will contribute one quadruple point
to $G_t$.
For type (2), any sheet of $D$ carries with it two sheets of 
$F'$, and so each sheet of $D$ involved in a quadruple point of $H'_t$ 
will double the number of 
corresponding quadruple points counted for $G_t$.
(e.g. if all four sheets involved in a quadruple point of $H'_t$ come from
$D$, then this occurrence will contribute $2^4 = 16$ quadruple points
to $G_t$.) And so there is an even number of quadruple points for $G_t$
in addition to those inherited from $H_t$.
And so $q(H_t)=q(G_t)=0$.
\end{pf}

The proof of the following lemma is a gathering of arguments from [W]:

\begin{lemma}\label{wh}
Let $M$ be an irreducible orientable 3-manifold, let $F$ be a (connected) 
closed 
orientable surface and $H_t:F\to M$ a
homotopy such that $H_0$ and $H_1$ are 
two incompressible embeddings with $H_0(F)=H_1(F)$.
Then there is a 
generic regular homotopy $G_t:F\to M$ 
such that 
the maps $H,G:\Q\to M$ defined by $H_t$ and $G_t$ are homotopic relative
$\pa (\Q)$, and such that the highest multiplicities of $G_t$ are 
double curves.
\end{lemma}

\begin{pf}
Denote $F' = H_0(F)$ ($=H_1(F)$.)
We may homotope $H$ relative $\pa (\Q)$ so that it will be transverse 
with respect to $F'$
and $H^{-1}(F')$ will be a system of incompressible surfaces in $\Q$.
By further homotoping we may assume 
$H^{-1}(F') = F\x\{s_0,s_1,...,s_n\}$ where $0=s_0<s_1<...<s_n=1$.

Let $p\in F$ be a chosen base point for $F$  and 
let $p'=H_0(p)$ be chosen as basepoint for $M$. 
We will think of $\p(F)=\p(F,p)$ as contained in $\p (M)=\p(M,p')$
via ${H_0}_*$.
Let $k:\hat{M}\to M$ be the covering corresponding to $\p(F)$
and let 
$\hat{F}=k^{-1}(F') \su \hat{M}$.
Let $\hat{H}:\Q\to \hat{M}$ be the lifting of $H$, then 
$\hat{H}^{-1}(\hat{F})=H^{-1}(F')$.
Let $\f_i$ be the component of $\f$ with $\hat{H}(F\x s_i)\su \f_i$.

We proceed by induction on $n$. Let $n=1$ and so $s_0=0, s_1=1$.
We distinguish two cases:

1. $\f_0=\f_1 :$ $\f_0=\hat{H}(F\x 0)$ is a strong deformation retract 
of $\hat{M}$. Let $J_s:\hat{M}\to\hat{M}$ 
denote this defomation, then $k\circ J_s\circ \hat{H}$
homotopes $H$ into $F'$
relative $\pa (F\x [0,1])$. And so now $H$ is a homotopy in
$F'$ between
$H_0$ 
and $H_1$. 
Such a homotopy is homotopic relative $\pa (\Q)$ to an isotopy, and we are 
done.

2. $\f_0\neq\f_1 :$ Cut $\hat{M}$ along $\f$ and let 
$\hat{M}'$ be the piece containing $\hat{H}(\Q)$.
By Lemma 5.1 of [W]
$\hat{M}' = \Q$.
And so 
we may homotope $\hat{H}$ relative
$\pa (\Q)$ so that $\hat{H}:\Q\to\hat{M}'$ 
will be a homeomorphism. 

Cut $M$ along $F'$ and let $M'$ be the piece covered by
$\hat{M}'$. (Perhaps there is no other piece.)
Then $H=k\circ\hat{H}:\Q\to M'$ is
now a covering map. 
Since $F'$ was originally covered exactly twice by $\pa (\Q)$, 
$H:\Q\to M'$ is a covering of degree at most 2. And so $H$, when thought of
again as a homotopy, is a regular homotopy which
has at most double points. 

We now let $n>1$ 
and so $0<s_1<1$. If now $\f_0=\f_1$ then 
(since $\f_0$ is a strong deformation retract of $\hat{M}$,)
we may homotope $\hat{H}$ relative
$\pa (\Q)$ as to push $F\x [0,s_1]$ to the other side of $\f_0$, 
by that reducing $n$. And so we may assume
$\f_0\neq\f_1$. As before we may use Lemma 5.1 of [W], to see that 
$\f_1$ is parallel to $\f_0$, and so we may now homotope $\hat{H}$ relative
$\pa (\Q)$ and without changing $\hat{H}^{-1}(\hat{F})$ such that
$\hat{H}|_{F\x s_1}$ will be an embedding, and then that 
$\hat{H}|_{F\x [0,s_1]}$ 
will be an embedding. 

We will now show that also
$H|_{F\x s_1} = k\circ \hat{H}|_{F\x s_1}$ is an embedding
(which is trivial if $F$ is not a torus.)
As before, $H:F\x [0,s_1]\to M'$ is a covering map. 
If in $M'$,
$H(F\x 0) \neq H(F\x s_1)$ then since $H|_{F\x 0}$ is an embedding,
$H(F\x 0)$ is covered exactly once, and so the same must be with 
$H(F\x s_1)$  and so $H|_{F\x s_1}$ is an embedding. If 
$H(F\x 0) = H(F\x s_1)$ and $H|_{F\x s_1}$ is some non-trivial covering,
then we will get via $H|_{F\x [0,s_1]}$ 
that $\p(F)$ is conjugate in $\p(M')$ to a 
proper subgroup of itself, which is impossible since it has finite index
(since it corresponds to the covering space $F\x [0,s_1]$.)

And so
$H|_{F\x s_1}$ is an embedding, and we may apply the induction hypothesis 
to $H|_{F\x [0,s_1]}$ and $H|_{F\x [s_1,1]}$.

\end{pf}

We are now ready to prove:

\begin{thm}\label{t3}

Let $M$ be an orientable irreducible 3-manifold with $\ppp (M) =0$.
Let $F$ be a system of closed orientable surfaces. 
If $H_t:F\to M$ is any CGRH in the regular homotopy class of an
embedding,
then $q(H_t)=0$.

\end{thm}

\begin{pf}
By the opening remark in
the proof of Proposition \protect\ref{p1}
and by induction on Lemma \protect\ref{cd}, 
we may assume that $H_0$ is an embedding and that for each
$i$ either $F^i=S^2$ or $H_0|_{F^i}$ is an incompressible embedding
(where $F^1, ... , F^n$ are the connected components of $F$.)

Given $p^i\in F^i$ we will 
construct a CGRH $G^i_t:F^i\to M$ satisfying: 
(1) $G^i_0=H_0|_{F^i}$.
(2) The paths 
$t\mapsto H_t(p^i)$ and $t\mapsto G^i_t(p^i)$ are homotopic 
relative endpoints. 
(3) The highest multiplicities of $G^i_t$ are double curves.
We then define 
$G_t:F\to M$ to be the CGRH which performs each $G^i_t$ in it's turn, 
while fixing all other components. By condition (3) the highest
multiplicities of $G_t$ will be triple points, and so $q(G_t)=0$.
By conditions (2), Proposition \protect\ref{p1} applies to
$H_t * G_{-t}$ and so $q(H_t)=q(H_t * G_{-t})=0$.

And so we are left with constructing the CGRH's $G^i_t$ satisfying 
(1), (2) and (3) above.

For $F^i=S^2$, $G^i_t$ is defined as follows: 
$H_0(S^2)$ bounds a ball $B$. Shrink $S^2$ inside $B$ 
until it is very small. Then move this tiny sphere along the loop
$t\mapsto H_t(p^i)$, finally re-entering $B$. Since $M$ is orientable,
our $S^2$ will return with the same orientation it originally had 
(One may think of a little ball 
bounded by $S^2$ which is moving along with it.)
And so we may 
isotope $S^2$ in $B$, back to its original position. And so $G^i_t$
is actually an isotopy.

For $H_0(F^i)$  incompressible, 
use Lemma \protect\ref{wh}.

\end{pf}

\begin{remark}

In the proof of Theorem \protect\ref{t3}, when we were assuming $F^i=S^2$ we
didn't actually need $H_0(S^2)$ to bound a ball. It would be enough that 
$\pp(M)=0$ since then we could deform $S^2$ by some generic regular 
homotopy $G_t$ until it is 
contained in a ball, inside which $G_t$ continues until the sphere is
embedded again.
Then move the embedded sphere around the loop, 
and then return to the original position by $G_{-t}$.
This would still contribute $q=0$.

And so, if the embedding $f:F\to M$ is such that repeated compression of
$f(F)$ will turn it into a union of spheres, 
then the theorem will still be true with the
(perhaps) weaker assumption $\pp(M)=0$ in place of irreducibility.
In particular, 
if $M$ is a fake open 3-cell, then 
the compressibility down to spheres always holds.
(On the other hand, when applying Lemma \protect\ref{wh}
to the case $H_0(F^i)$ is incompressible,
we indeed use irreducibility.)
\end{remark}

\section{Surfaces in 3-Manifolds: Counter-Examples}

1. The most obvious example of a CGRH with $q=1$ is 
the following $H_t:S^2\to S^3$. $H_0$ 
is some embedding. $H_t$ then isotopes 
$S^2$ in $S^3$ as to reverse the sides of $S^2$. Finally
it performs a generic
eversion inside a ball 
$B\su S^3$ so as to return $S^2$ to its original position.
By Theorem \protect\ref{mb2} $q(H_t)=1$. 
Furthermore, this $H_t$ may be performed while fixing some disc in 
$S^2$. And so $H_t$ satisfies all the conditions of
Theorems \protect\ref{t2} and \protect\ref{t3} except $\ppp(M)=0$.

2. Let $h:S^3\to M$ be the covering map where $M=\R P^3$ or $L_{3,1}$,
and so $h$ is a double or triple covering. 
Perform the above $H_t:S^2\to S^3$ such that a ball $B$ in which
we start $H_t$ and in which
the eversion takes place, is embedded into $M$ by $h$. 
Let $G_t=h\circ H_t$.
Since the covering is double or triple, the first part of 
$H_t$ contributes no quadruple points, and so
$q(G_t)=1$.

3. Lemma \protect\ref{l2}
provides an example of $q=1$ if one takes $F$ to be a surface of
odd Euler characteristic. Instead of having $D$ fixed, we let it rotate once
so as to cancel the Dehn twist. And so we get a CGRH with 
$q=\chi(F) \ mod \ 2 \ = 1$. The condition of Theorem
\protect\ref{t2} which is violated is that no disc is fixed, 
though, we may perform this CGRH while keeping just a point fixed 
(some point in $D$ around which we rotate $D$ to cancel the Dehn twist.) 

4. In any non-orientable 3-manifold there is a CGRH with $q=1$. Take a little
sphere inside a ball. Move the little sphere along an orientation reversing 
loop until it returns to itself with opposite orientation. 
Then perform an eversion inside the ball to return to the original position.

5. Let $F_1$, $F_2$ be two surfaces such that there exists an immersion 
$g:F_2\to F_1\x S^1$ with an odd number of triple points.
(e.g. such an immersion exists whenever $F_2$ has odd Euler characteristic.)
Let $F=F_1\cup F_2$ and $M=F_1 \x S^1$. Let $H_0:F\to M$ be $g$ on
$F_2$ and the inclusion $F_1\to F_1\x *\su F_1\x S^1$ on $F_1$.
Let $H_t$ move $F_1$ once around $M$ while fixing $F_2$, then
$q(H_t)=1$.

\section{Embeddings of Tori in $\E$}

Let $T$ be the standard torus in $\E$.
For the sake of definiteness say $T$ is the surface obtained
by rotating the circle $\{ \ y=0 \ , \ (x-2)^2+z^2=1 \ \}$ about the 
$z$-axis. Let $m\su T$ be the 
meridian, e.g. the loop $\{ \ y=0 \ , \ (x-2)^2+z^2=1 \ \}$. 
$m$ bounds a disc in the compact side of $T$.
Let  $l\su T$ be the longitude, e.g. the loop 
$\{ \ z=0 \ , \ x^2+y^2=1 \ \}$. 
$l$ bounds a disc in the non-compact side of
$T$. We choose $m,l$ (with some orientation)
as the basis of $H_1(T,\Z)$, and this induces an
identification between $H_1(T,\Z)$ and $\Z^2$. The group $M(T)$
of self diffeomorphisms of $T$ up to homotopy, is then identified with
$GL_2(\Z)$, and so $f\in M(T)$ will be thought of both as a map
$T\to T$ and as a $2\x 2$ matrix. Let $\tau$ denote reduction mod 2,
so we have $\tau:\Z^2\to(\C)^2$ and $\tau:GL_2(\Z) \to GL_2(\C)$.
We will use round brackets for objects over $\Z$, and square brackets
for objects over $\C$.
Let $H\su GL_2(\C)$ be the subgroup $\{ \ \U \ , \ \V \ \}$.

\begin{prop}\label{mcg}
Let $T\su \E$ be the standard torus, and denote the inclusion map
by $i:T\to\E$. Let $M(T)$ be identified with $GL_2(\Z)$ 
via $m,l$ as above and let
$H\su GL_2(\C)$ be as above. 

(1) For $f\in M(T)$,  $i\circ f$
is regularly homotopic to $i$ iff $\tau(f)\in H$.

(2) If $\tau(f)=\U$ then any generic regular homotopy $H_t:T\to\E$
between $i\circ f$ and $i$ will have $q(H_t)=0$.
If $\tau(f)=\V$ then any generic regular homotopy $H_t:T\to\E$
between $i\circ f$ and $i$ will have $q(H_t)=1$.

\end{prop}

\begin{pf}
Let $V_1,V_2$ be the standard 2-frame on $T$, i.e. $V_1$ (resp. $V_2$)
is everywhere
tangent to the translates of $m$ (resp. $l$.)
Let $K(v_1,v_2)$ be the function 
defined in the beginning of the proof
of Proposition \protect\ref{t0}. For an immersion $g:T\to \E$, let 
$h_g:T\to SO_3$ be defined by
$h_g = K(dg(V_1) , dg(V_2))$.
By the Smale-Hirsch Theorem, 
$i\circ f$ will be regularly homotopic to $i$ iff
$h_{i\circ f}$ is homotopic to $h_i$.
It is easy to see that
$h_{i\circ f}$ is homotopic to $h_i\circ f$, and so
$i\circ f$ will be regularly homotopic to $i$ iff
$h_i\circ f$ is homotopic to $h_i$.
Let $k:H_1(T,\Z)\to H_1 (SO_3,\Z)=\C$ be the 
homomorphism induced by $h_i$. Since $\pp(SO_3)=0$, 
$h_i\circ f$ will be homotopic to $h_i$ iff 
$k\circ f = k$. (Recall that we denote by $f$ both the map $T\to T$
and the automorphism it induces on $H_1(T,\Z)$.)
Let $k':H_1(T,\C)\to \C$ be defined by $k=k'\circ\tau$,
then $k\circ f = k$ iff $k'\circ \tau(f) = k'$. 
It is easy to see that $k{{\left(
\begin{smallmatrix}
1 \\
0
\end{smallmatrix}
\right) }}
= k{{\left(
\begin{smallmatrix}
0 \\
1
\end{smallmatrix}
\right) }}=1$
and so $k'\e=k'\ee=1$ and $k'\eee = 0$. 
And so $k'\circ \tau(f) = k'$ iff $\tau(f)$ maps the set
$\{ \ \e , \ee \ \}$ into itself.
Assertion (1) follows.

It is easy to verify, by means of row and column operations, that the matrices
$\left(
\begin{smallmatrix}
1 & 2 \\
0 & 1
\end{smallmatrix}
\right)$,
$\left(
\begin{smallmatrix}
1 & 0 \\
2 & 1
\end{smallmatrix}
\right)$,
$\left(
\begin{smallmatrix}
-1 & 0 \\
0 & -1
\end{smallmatrix}
\right)$
and
$\left(
\begin{smallmatrix}
-1 & 0 \\
0 & 1
\end{smallmatrix}
\right)$
generate $\tau^{-1}\U$.
By Theorem \protect\ref{t1}, in order to prove the first half of 
assertion (2), we only need to construct one regular homotopy for
each one of the four generators, and see that it has $q=0$.

$\left(
\begin{smallmatrix}
1 & 2 \\
0 & 1
\end{smallmatrix}
\right)$ is a double Dehn twist along the meridian. 
The regular homotopy will be as follows: Perform two $A$ moves along
parallel meridians, each giving a Dehn twist of the required orientation,
and such that the rings are formed on the same side. Then perform 
a $B$ move to merge the two rings. We now have a ring bounding a disc
in $T$. We may get rid of it with one 
$A^{-1}$ move. And so we are left with 
precisely the two Dehn twists we needed.  As to $q$, we had two $A$ moves, 
one $B$ move and one $A^{-1}$ move and so by Lemma \protect\ref{l1}
we indeed have $q=0$. (Actually, if we would have formed the rings on 
opposite sides of $T$ then we could have completed the regular 
homotopy with no additional quadruple points. So we would have just
the two $A$ moves which would indeed also give $q=0$.)

$\left(
\begin{smallmatrix}
1 & 0 \\
2 & 1
\end{smallmatrix}
\right)$ is a double Dehn twist along
the longitude. Since 
for the standard torus, the longitude is also disc 
bounding (on the non-compact side,) we can proceed exactly as in the 
previous case.

$\left(
\begin{smallmatrix}
-1 & 0 \\
0 & -1
\end{smallmatrix}
\right)$ may be achieved by rigidly rotating $T$ around the $x$ axis through 
an angle of
$\pi$, and so again $q=0$.

$\left(
\begin{smallmatrix}
-1 & 0 \\
0 & 1
\end{smallmatrix}
\right)$ is a reflection with respect to the $xy$-plane. We can achieve 
it by a regular homotopy as follows: Perform two $A$ moves on two longitudes,
say the circles $\{ \ z=0 \ , \ x^2+y^2=1 \ \}$
and $\{ \ z=0 \ , \ x^2+y^2=3 \ \}$, with the rings both facing the 
non-compact side, and with
the orientation
of the Dehn twists chosen so that they cancel each other.
We may then continue with just double curves, to exchange the upper
and lower halves of $T$ until we arrive at the required reflection.
We had two $A$ moves and so $q=0$.

And so we have proved the first half of assertion (2). Since 
$\tau^{-1}\U$ has index 2 in $\tau^{-1}(H)$, we need in order to prove the 
second half of (2), just to check for one element. We construct the following
regular homotopy. Start with an $A$ move on some small circle that
bounds a disc in $T$, and with the ring facing the 
non-compact side. Let $u$ be a circle on the ring, which is parallel 
to the intersection circle. It bounds a disc $U$ of $T$. 
Keeping $u$ fixed, push and expand $U$ all the way to the other side of $T$
so that it encloses all of $T$, so the torus is now embedded again.
(This last move required only 
double curves.) A disc that spanned $m$ is now on the 
non-compact side, and a disc that spanned $l$ is now
on the compact side. And so if we now isotope $T$ until it's image 
again coincides with itself, we will have an $f$ which interchanged
$m$ and $l$, and so $\tau(f)=\V$. This regular homotopy
required one $A$ move, then some double curves and finally an isotopy,
and so $q=1$.
\end{pf}

We continue to think of $T$ as the standard torus contained 
in $\E$, with inclusion
map $i$.
Now
let $f:T\to \E$ be any embedding,
then $f(T)$ again separates $\E$ into 
two pieces, one compact and the other non-compact.
Denote the compact
piece by $C_f$, and the non-compact piece by $N_f$.

For a moment, compactify $\E$ with a point $\infty$, so that
$\E\cup \{\infty\}=S^3$.
$f(T)\su S^3$ always bounds a solid torus on (at least) one side.
The other side is then a knot complement. There are now two cases: 
(a) $\infty$ is in the knot complement side, and so $C_f$ is a solid torus 
(Figure 3a.) (b) $\infty$ is in the solid torus side, and so $C_f$ is 
a knot complement (Figure 3b.) It is of course possible to have a solid torus 
on both sides, in which case the knot in question is the unknot
(Figure 3c.)

\begin{figure}[h]
\ \ \ \ \ \ \ \ 
\scalebox{.6}{\includegraphics{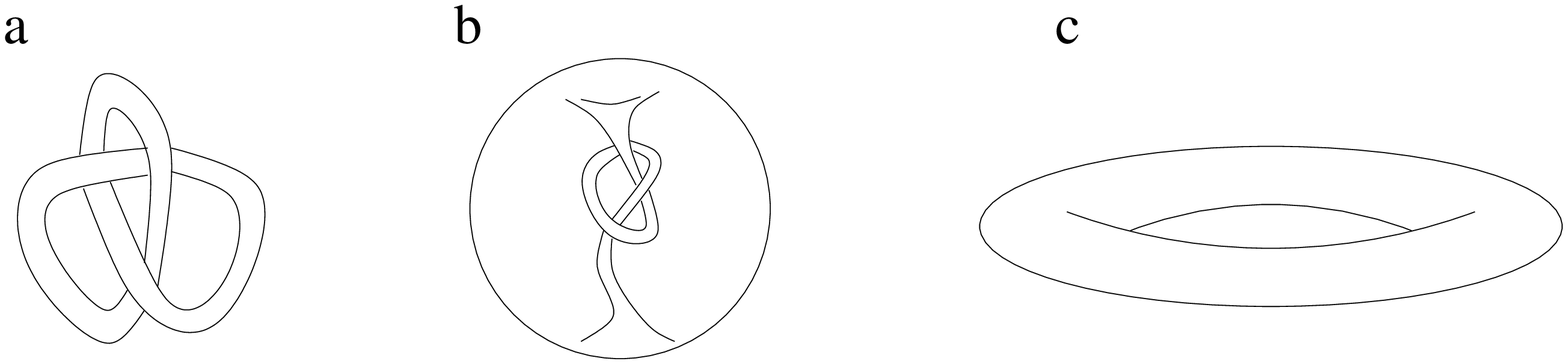}}
\caption{Embeddings of $T$ in $\E$}
\end{figure}

$f$ maps $T$ into both $C_f$ and $N_f$ and so we have maps
$H_1(T,\C)\to H_1(C_f,\C)$ and $H_1(T,\C)\to H_1(N_f,\C)$.
Denote the unique non-zero element of the kernel of each of these maps
by $c_f$ and $n_f$ respectively. By the 
above observation that $f(T)$ may be
thought of as the boundary of a regular neighborhood of some
knot $k\su S^3$, we see that $c_f$ and $n_f$ generate $H_1(T,\C)$, 
which means simply that they are two distinct 
elements of the set $\{ \ \e , \ee , \eee \ \}$.
Let us now choose once and for all some arbitrary ordering of this 
set, say $\e < \ee < \eee$. 
Define $Q(f)=0 \in \C$ 
if $c_f < n_f$ and $Q(f)=1 \in \C$ if $c_f > n_f$.
Now let 
$F=T^1 \cup ... \cup T^n$ be a union of 
$n$ copies of $T$.
If $f:F\to \E$ is an embedding, then define
$Q(f)= \sum_{i=1}^n Q(f|_{T^i}) \in \C$.

We will now prove:

\begin{thm}\label{tori}
Let $F$ and $Q$ be as above, 
and let $f,g:F\to \E$ be two embeddings.
If $f$ and $g$ are regularly homotopic,
then any generic regular homotopy $H_t$ between them will satisfy
$q(H_t) = Q(f) - Q(g)$.
\end{thm}

\begin{pf}

By Theorem \protect\ref{t1} any such $H_t$ will have the same $q$.
And so we need to verify the theorem for just one $H_t$.

Since $f$ is an embedding, if we rigidly move one of the $T^i$s 
while keeping the others fixed, 
only double curves will appear. And so we may
move them one by one until they are contained in disjoint balls,
and this will contribute nothing to $q$. 
Since $Q$ is defined by computing it
for each component separately, $Q(f)$ is also unchanged.
We may do the same with $g$, and so we actually need to deal with 
each component separately, and so we may assume from now on that $F$ is 
the one torus $T$. 

As mentioned, $f(T)$ is either (a) a torus which bounds a solid torus,
or (b) a sphere with a tube running inside it. 
See Figure 3a and 3b respectively. 
We now begin our regular homotopy by having the solid torus, in case (a),
or the inner tube, in case (b), pass itself until we eliminate all knotting. 
This creates only double curves and so
contributes nothing to $q$, and we claim that it also does not change 
$c_f$ and $n_f$. This follows from the fact that the $\C$-meridian and 
the $\C$-longitude of a knot, do not change
under such crossing moves. (In case (a) $c_f$ is the $\C$-meridian
of a knot and $n_f$ is the $\C$-longitude, and in case (b)
it is the other way around.) And so since $c_f$ and $n_f$ remain unchanged,
$Q$ is unchanged. We may now continue with an isotopy
until the image of $T$ coincides with $T$ itself.

As we may do the same for $g$, we may assume from now on that 
the image of both $f$ and $g$ is $T$ itself,
and so $f=i\circ f'$ and $g=i\circ g'$ for some $f',g':T\to T$.
If $H_t$ is a regular homotopy from $f$ to $g$, then
$H_t \circ {f'}^{-1}$ is a regular homotopy from $i$ to 
$i\circ (g' \circ {f'}^{-1})$.
By Proposition 
\protect\ref{mcg}(1) $\tau (g'\circ {f'}^{-1})$ is either $\U$ or
$\V$. By \protect\ref{mcg}(2),
if $\tau(g'\circ {f'}^{-1}) = \U$ (i.e. $\tau(g')=\tau(f')$,) then
$q(H_t \circ {f'}^{-1}) =0$ 
and if $\tau(g'\circ {f'}^{-1}) = \V$ (i.e. 
$\tau (g')=\V\circ\tau(f')$) then $q(H_t \circ {f'}^{-1})=1$.
 
If $\tau(f')=\tau(g')$ then $c_f=c_g$ and $n_f=n_g$ and so 
$Q(f)=Q(g)$ and so $Q(f)-Q(g)=0=q(H_t\circ {f'}^{-1}) = q(H_t)$
and we are done.

Assume now $\tau (g')=\V\circ\tau(f')$. 
By definition of $c_f$ and $n_f$ we 
must
have in $\C$ homology: $\tau(f')(c_f)=\tau(g')(c_g)=[m]$ and 
$\tau(f')(n_f)=\tau(g')(n_g)=[l]$. 
So  $\tau(g')(c_f) = \V\circ\tau(f')(c_f)=\V[m]=[l]=\tau(g')(n_g)$,
and so $c_f = n_g$. Similarly
$\tau(g')(n_f) = \V\circ\tau(f')(n_f)=\V[l]=[m]=\tau(g')(c_g)$
and so $n_f = c_g$. That is, the pairs $c_f,n_f$ and $c_g,n_g$ are the 
same pair of elements of $H_1(T,\C)$ just with opposite order.
So the order of one pair matches the chosen order iff
the order of the other pair doesn't.
And so $Q(f)-Q(g)=1=q(H_t \circ {f'}^{-1})=q(H_t)$.
\end{pf}

\end{document}